# POLLING SYSTEMS WITH PARAMETER REGENERATION, THE GENERAL CASE

By Iain MacPhee, Mikhail Menshikov,[1]
Dimitri Petritis[2] and Serguei Popov[3]

*University of Durham, University of Durham, Université de Rennes 1 and Universidade de São Paulo*

We consider a polling model with multiple stations, each with Poisson arrivals and a queue of infinite capacity. The service regime is exhaustive and there is Jacksonian feedback of served customers. What is new here is that when the server comes to a station it chooses the service rate and the feedback parameters at random; these remain valid during the whole stay of the server at that station. We give criteria for recurrence, transience and existence of the $s$th moment of the return time to the empty state for this model. This paper generalizes the model, when only two stations accept arriving jobs, which was considered in [*Ann. Appl. Probab.* **17** (2007) 1447–1473]. Our results are stated in terms of Lyapunov exponents for random matrices. From the recurrence criteria it can be seen that the polling model with parameter regeneration can exhibit the unusual phenomenon of null recurrence over a thick region of parameter space.

**1. Introduction.** There is a large literature, for example, [1, 7, 9, 12, 14, 15] on polling systems, that is, systems with a single moving server and multiple stations where queues form. Systems of this kind arise frequently in production or computer networks (e.g., consider the situation when a single processing unit has to deal with job requests coming from several terminals). We consider a polling system where the queues at stations can be arbitrarily long and the service regime is exhaustive, that is, the server does not leave

Received March 2007; revised October 2007.
[1]Supported in part by FAPESP (2004/13610–2).
[2]Supported in part by the "Réseau Mathématique France–Brésil" and by the European Science Foundation.
[3]Supported in part by CNPq (302981/02–0), FAPESP (04/03056–8 and 04/07276–2) and by the "Rede Matemática Brasil–França."

*AMS 2000 subject classifications.* Primary 60K25, 60J10; secondary 60G42, 90B22.
*Key words and phrases.* Polling system, parameter regeneration, stability, time-inhomogeneous Markov chains, recurrence, Lyapunov functions, random matrices.

a station until the queue there is empty. Also, customers can move from one queue to another one before definitively leaving the system. What is new here and in [11] is that we allow a randomly selected service regime each time the server visits a station. Throughout the paper we use *service regime* to mean the set of service and feedback parameters currently in use; in other words, we are working with models where those parameters are (randomly) updated each time the server goes to another station. It is worth mentioning that normally in the polling literature arrival and service parameters are fixed, and our main purpose here is to explore new phenomena which appear when one allows them to be randomly updated. We mention also the related work [3], where a single server queue with two possible regimes was considered. One of the results of [3] is that the asymptotic behavior of the process is very different when the slow service rate is less than the arrival rate of jobs than when it is faster than the arrival rate. This seems to be very natural from our point of view; cf. Theorem 1.1 and Corollary 1.1.

This paper presents results for a generalization of the model studied in [11]. The model considered there is essentially two dimensional, while here we consider any number of stations. The reason for splitting this work into two separate papers, depending on whether two or more stations receive jobs, is that in the former case the stability conditions are expressible in terms of *explicitly* determined factors, while in the latter, they are expressible in terms of factors whose existence is shown but whose values are inaccessible to direct computation. These factors are reminiscent of the Lyapunov exponents for products of random matrices whose values usually can be approximated only by stochastic simulation.

The model in [11] includes general service time and switching time distributions. With more than two stations the model and its analysis become considerably more complex, so in this paper we discuss exponential service times, zero switching times and a cyclic server routing scheme to help simplify the discussion. Within this setup, the results of the present paper generalize those of [11] and the proofs differ substantially in detail as well.

It seems to us that queueing models in which parameters are randomly reset are a very natural generalization of standard queueing models. However, stability investigations for such models seem to require quite deep mathematical techniques as these models can exhibit new phenomena. In particular, our polling model can have a "thick" region of null recurrence and even in ergodic cases the return time to the empty state may have only polynomial moments. This implies that convergence to equilibrium cannot occur at exponential rate.

1.1. *The model.* Consider the following continuous time model: there are $d+1$ nodes (stations) indexed by $n \in \mathbb{Z}_{d+1} := \{0, \ldots, d\}$, each one has its own infinite capacity queue. Customers arrive at the queue of station $n$ according to a Poisson point process with rate $\lambda_n > 0$, independently of all other events, $n = 0, \ldots, d$. A single server visits the stations in cyclic order $(0 \to 1 \to \cdots \to d \to 0 \to \cdots)$ starting at station 0. At each station the server adopts the *exhaustive* serving policy: it serves all the customers that were queueing at the station when the server arrived there together with all subsequent arrivals up until the queue first becomes empty and then instantly jumps to the next station. Every served customer either leaves the system or is transferred to another queue, with given probabilities. Each time the server arrives at a station, it chooses randomly a serving rate and feedback (transition) probabilities; those parameters remain valid during the stay of the server at that station.

To make the last sentence rigorous, we have to make some definitions. Abbreviating $(\mu, \gamma_1, \ldots, \gamma_d)$ by $(\mu, \gamma)$, define

$$\mathcal{M} = \left\{ (\mu, \gamma) : \mu > 0, \gamma_i \geq 0 \text{ for all } i, \text{ and } \sum_{j=1}^{d} \gamma_j \leq 1 \right\} \subset \mathbb{R}^{d+1}.$$

Let $\nu_n$, $n = 0, \ldots, d$, be $d+1$ probability measures on $\mathcal{M}$. Denote also $[k] := k \bmod (d+1) \in \mathbb{Z}_{d+1}$. Then, we complete the description of the model as follows: every time the server comes to station $n$ we choose a random element $(\mu, \gamma) \in \mathcal{M}$ according to $\nu_n$, independently of everything else. For the period while the server performs this batch of services, the service times are i.i.d. and exponentially distributed with parameter $\mu$. Each served customer either goes to station $[n+k]$ with probability $\gamma_k$ for $k = 0, \ldots, d$ or leaves the system with probability $1 - \gamma_1 - \cdots - \gamma_d$, independently of other events. Also, for definiteness we assume that the server just waits somewhere (say, at some special place $\mathfrak{S}$) when there are no customers in the system (but when the first customer arrives, the server instantaneously jumps to the corresponding station).

Throughout this paper, and often without recalling it explicitly, we suppose that the (uniform ellipticity) condition below is fulfilled:

CONDITION E.

(i) There exists $\varepsilon_0 > 0$ such that, for all $n = 0, \ldots, d$,
$$\nu_n[(\mu, \gamma) : \mu > \lambda_n + \varepsilon_0] = 1.$$

(ii) There exists $M_0$ such that, for all $n = 0, \ldots, d$,
$$\nu_n[(\mu, \gamma) : \mu \leq M_0] = 1.$$

Part (i) of the above condition guarantees that the server cannot get stuck in a single station. Part (ii) is just a technical assumption [using a simple monotonicity argument, one can show that this assumption is not necessary for proving parts (ii) and (iii) of Theorem 1.1 below].

1.2. *The main result.* Assuming that the process starts from some fixed and nonempty initial configuration (and initially the server is at station 0), we define $\tau$ as the first moment when there are no customers in the system (i.e., $\tau$ is the time of reaching the zero configuration).

To formulate our result, we need to introduce certain random matrices. For $n \in \{0, \ldots, d\}$, let $A^{(n)} = (a_{i,j}^{(n)})_{i,j=0,\ldots,d-1}$ be a matrix defined as follows:

$$a_{0,k}^{(n)} = \frac{\lambda_{[n+k+1]} + \gamma_{k+1}\mu}{\mu - \lambda_n},$$

$k = 0, \ldots, d-1$, $a_{k,k+1}^{(n)} = 1$, $k = 0, \ldots, d-2$, and $a_{i,j}^{(n)} = 0$ for all other $i, j$, where $(\mu, \gamma) \sim \nu_n$. That is,

$$(1.1) \qquad A^{(n)} = \begin{pmatrix} \frac{\lambda_{[n+1]} + \gamma_1 \mu}{\mu - \lambda_n} & 1 & 0 & 0 & \ldots & 0 \\ \frac{\lambda_{[n+2]} + \gamma_2 \mu}{\mu - \lambda_n} & 0 & 1 & 0 & \ldots & 0 \\ \frac{\lambda_{[n+3]} + \gamma_3 \mu}{\mu - \lambda_n} & 0 & 0 & 1 & \ldots & 0 \\ \vdots & \vdots & \vdots & \vdots & \ddots & \vdots \\ \frac{\lambda_{[n+d-1]} + \gamma_{d-1}\mu}{\mu - \lambda_n} & 0 & 0 & 0 & \ldots & 1 \\ \frac{\lambda_{[n+d]} + \gamma_d \mu}{\mu - \lambda_n} & 0 & 0 & 0 & \ldots & 0 \end{pmatrix}.$$

Let $A = A^{(d)} \ldots A^{(0)}$, where $A^{(0)}, \ldots, A^{(d)}$ defined above are *independent* random matrices. Throughout the paper, we suppose that the random matrix $A$ is really random, that is, one cannot find a deterministic matrix $B$ such that $A = B$ a.s. (this means, in fact, that in at least one node there is randomness with respect to the choice of the regime).

REMARKS.

- The motivation for the definitions of the random matrices $A^{(n)}$ and $A$ is given in the arguments leading to (2.9) below. Very briefly, for any given queue length vector $x$, the random vector $Ax$ gives the probable queue lengths after the server has completed a cycle of visits to all stations (a random vector because the service rates and feedback probabilities are random).
- This random matrix $A$ plays a role similar to the matrix $M(s)$ used in our previous paper [11], but there are important differences. $M(s)$ is only a matrix in [11] because we permit Markov switching of the server there. With the cyclic server route assumed in this paper, $M(s)$ would reduce to a single random variable, which is why we are able to state our criteria there in terms of deterministic parameters.

In the following few lines we introduce the main classification parameter. Let $A_1, A_2, A_3, \ldots$ be a sequence of independent random matrices, each one having the same distribution as $A$. Define for any $s \geq 0$

$$\mathsf{k}(s) = \lim_{n \to \infty} (\mathbb{E}\|A_n \ldots A_1\|^s)^{1/n} \tag{1.2}$$

(unless otherwise stated, all norms we use in this paper are $\mathcal{L}_1$-norms). As shown below, this limit always exists by submultiplicativity of the norm. Indeed, let $I = \{s \in \mathbb{R} : \mathbb{E}\|A\|^s < \infty\}$. Observe that Condition 1.1 guarantees that $[0, \infty) \subseteq I$.

Since the map $A \mapsto \|A\|$ is submultiplicative, $\mathsf{k}_n(s) = \mathbb{E}\|A_n \ldots A_1\|^s$ satisfies $\mathsf{k}_{n+m}(s) \leq \mathsf{k}_n(s)\mathsf{k}_m(s)$, hence, $\lim_n \frac{\log \mathsf{k}_n(s)}{n} = \log \mathsf{k}(s)$ exists on $I$.

If $s_1, s_2 \in I$ with $s_1 < s_2$, then the whole segment $[s_1, s_2] \subseteq I$. For any $t \in (0, 1)$, by the Hölder inequality with $p = 1/t$ and $q = 1/(1-t)$, we have

$$(\mathsf{k}_n(ts_1 + (1-t)s_2))^n \leq \mathbb{E}(\|A_n \ldots A_1\|^{pts_1})^{1/p} \mathbb{E}(\|A_n \ldots A_1\|^{q(1-t)s_2})^{1/q}.$$

This inequality establishes log-convexity of $\mathsf{k}(s)$ [and shows also that $\mathsf{k}(s)$ is continuous]. Also, from the nondegeneracy of the random matrices, it follows that $\mathsf{k}(s)$ is nonconstant. Furthermore, $\mathsf{k}'(0)$ is the top Lyapunov exponent of $A$.

The main classification parameter is

$$s_0 = \inf\{s > 0 : \mathsf{k}(s) > 1\} \tag{1.3}$$

(by definition, $\inf \varnothing = +\infty$). Note that $\mathsf{k}(0) = 1$, and so $s_0 = 0$ if $\mathsf{k}'(0) > 0$, and $0 < s_0 \leq +\infty$ if $\mathsf{k}'(0) < 0$. The main result of this paper is the following.

THEOREM 1.1. (i) *If $\mathsf{k}'(0) > 0$, then $\mathbb{P}[\tau = \infty] > 0$, that is, the process is transient.*
  (ii) *If $\mathsf{k}'(0) < 0$, then $\mathbb{P}[\tau = \infty] = 0$, that is, the process is recurrent.*
  (iii) *Suppose that $0 < s < s_0$, then $\mathbb{E}\tau^s < \infty$.*
  (iv) *For $s > s_0$, we have $\mathbb{E}\tau^s = \infty$.*

We do not consider the case when $\mathsf{k}'(0) = 0$. There is some discussion of the difficulties for the case of two stations in [11].

The above theorem shows that, for the model of this paper, the typical situation is that geometric ergodicity is absent, and only polynomial moments of the return time exist (this corresponds to the case $s_0 < \infty$). We now formulate a simple sufficient condition for the finiteness of $s_0$. For any matrix $B$, let $\rho(B)$ be the spectral radius of $B$ [for the case of nonnegative matrices, $\rho(B)$ is the maximal eigenvalue]. We have

$$\mathbb{P}[\rho(A) > 1] > 0 \quad \implies \quad s_0 < \infty. \tag{1.4}$$

Indeed, the fact that $\mathbb{P}[\rho(A) > 1] > 0$ implies, by continuity, that there exist $\varepsilon, \delta > 0$ and a nonnegative vector $y$ such that $\|y\| = 1$ and $\mathbb{P}[Ay > (1+\varepsilon)y] > \delta$. Then

$$\mathbb{E}\|A_n \ldots A_1\|^s \geq \mathbb{E}\|A_n \ldots A_1 y\|^s$$
$$\geq (1+\varepsilon)^{sn}\delta^n,$$

and so we obtain from (1.2) and (1.3) that if $s$ is such that $(1+\varepsilon)^s \delta > 1$, then $s_0 < s$.

The relation (1.4) has a nice interpretation:

COROLLARY 1.1. *For the absence of geometric ergodicity (i.e., for $s_0$ to be finite) it is sufficient that there exists a set of* possible *parameters [i.e., $(\mu, \gamma)_n \in \operatorname{supp}\nu_n$, $n = 0, \ldots, d$] for which the* homogeneous *model is transient.*

Indeed, note that for a nonrandom (and nonnegative) matrix, the logarithm of the largest eigenvalue equals the top Lyapunov exponent, and so by Theorem 1.1(i), the model with nonrandom regimes is transient if $\rho(A) > 1$.

REMARK. As with the model of [11], the model of this paper can display an interesting feature, the existence of a thick null recurrence phase. Think about a family of models depending on one (or several) parameter(s). Suppose that this family is nicely defined, in such a way that $s_0$ is a continuous function of those parameters, and such that 1 is in the interior of the range of values of this function. Then, it is clear that, as we continuously change the values of the parameters, there are nondegenerate regions that correspond to ergodicity/null recurrence. This kind of behavior is quite rare; in particular, when the parameters of the polling system are fixed (say, $\lambda_i$ is the arrival rate to station $i$, $\mu_i$ is the service rate when the server is in station $i$, and there is no feedback), it is known that null recurrence is a critical phenomenon which can only occur when $\sum_i \lambda_i/\mu_i = 1$; see [7]. Only a few examples of models having a thick null recurrence phase are known (besides [11], see also [5]). For a concrete example which demonstrates this situation in the setup of this paper, we refer the reader to Section 2 of [11]. The example there has only two stations and so the calculations of the required exponents can be done explicitly. Note that, for the case of two stations, the models of this paper and of [11] (with exponential services and without switching time) are the same.

1.3. *Related problems and generalizations.* First, we remark that the model with any deterministic and periodic itinerary of the server can be studied in exactly the same way as the model of this paper. We decided to focus on the circular itinerary to avoid notational complications. There are other interesting routing policies that could be studied, for example, Markovian (after finishing the service in node $i$, the server goes to $j$ with probability $P_{ij}$), greedy (the server goes to the largest queue), local greedy (the server goes to the queue that is largest among some set of neighbors), etc. In this general setting, however, such routing policies cannot be treated directly by the methods of this paper, so a modification (or at least a nontrivial refinement) of our approach would be necessary. For the case when only two stations can receive customers, Markovian routing was considered in [11].

We believe that the results of this paper remain valid if we allow the incoming flow rates to be random as well, updating them to new values (independently from the past) at each server switching time [of course, one has to require that a condition analogous to Condition 1.1(ii) holds also for the $\lambda$'s]. This leads to a lot more complexity in our proofs but no new phenomena, so we have not formulated the model in this way. Further, it did not seem very natural to us that switching events should affect the incoming rates at the stations which are "far" from the server. A more natural model would update the incoming rates only at the stations which are "close" to the server (similarly to [11]), but this introduces yet more modeling complications.

There are two other natural ways to generalize the present model: first, consider service time with general distributions and second, consider the case of nonzero switching time (i.e., when the server empties the queue on a station, it takes some time traveling to the next one). With some additional assumptions, it is possible to prove that the analogue of Theorem 1.1 remains valid. Again, we decided not to describe the present model in this general form, because it does not lead to any new phenomena and complicates the proofs a lot.

In the remaining part of this section we describe two other models for which our method works.

*Revolver model.* This time, consider a model with stationary server and moving queues. There are $d+1$ locations marked $0, \ldots, d$ and customers arrive at location $n$ as a Poisson stream of rate $\lambda_n$. The server is always in location 0 and serves the queue using a regime $(\mu, \gamma)$ (where $\mu$ is the rate of service, and $\gamma_n$ is the probability that a served customer feeds back into queue $n$, $n = 1, \ldots, d$). When the queue at location 0 is emptied, all the customers in the location $i$ are instantly transferred to $i - 1$, $i = 1, \ldots, d$, and the server chooses a new regime $(\mu', \gamma')$ independently of the past and according to some probability measure $\tilde{\nu}$ on $\mathcal{M}$.

That is, we can visualize this system in the following way: the queues are located on a wheel, and when the server finishes the current batch of service, the wheel turns and delivers the next queue to the server. Define [compare with (1.1)]

$$A' = \begin{pmatrix} \frac{\lambda_1 + \gamma_1 \mu}{\mu - \lambda_0} & 1 & 0 & 0 & \cdots & 0 \\ \frac{\lambda_2 + \gamma_2 \mu}{\mu - \lambda_0} & 0 & 1 & 0 & \cdots & 0 \\ \frac{\lambda_3 + \gamma_3 \mu}{\mu - \lambda_0} & 0 & 0 & 1 & \cdots & 0 \\ \vdots & \vdots & \vdots & \vdots & \ddots & \vdots \\ \frac{\lambda_{d-1} + \gamma_{d-1} \mu}{\mu - \lambda_0} & 0 & 0 & 0 & \cdots & 1 \\ \frac{\lambda_d + \gamma_d \mu}{\mu - \lambda_0} & 0 & 0 & 0 & \cdots & 0 \end{pmatrix},$$

where $(\mu, \gamma) \sim \tilde{\nu}$. Then [with the quantities $\mathsf{k}(s)$ and $s_0$ defined as in (1.2) and (1.3)] Theorem 1.1 remains valid for the revolver model. This can be proved in a similar way to the proof of Theorem 1.1 below (in fact, it is easier because in the revolver model the server is always in the same position, and so it is not necessary to introduce $d+1$ random matrices with different distributions).

*Model with a gated service discipline.* This is a modification of the model defined in Section 1.1, described as follows. Instead of the exhaustive service discipline, we consider the so-called *gated service discipline* (see, e.g., Section 3.2 of [1]): the server only serves the customers that were in the station on its arrival, and then proceeds to the next station. The service regime can be described by the same set of parameters $(\mu, \gamma)$, and here even the inclusion of parameter $\gamma_0$ (which corresponds to the probability of sending the customer to the same station, but *before* the gate) makes sense. For this model, we still need an ellipticity condition, but in a less restrictive form. For gated service, we only need to require that $\varepsilon_0 < \mu < M_0$ with probability 1. This is because, when the server arrives at the station, it only needs to serve a fixed number of customers (those who were there at the moment of arrival), and so, even if it chooses a very small $\mu$, the service will be eventually finished, and the server will go to the next station (leaving behind the clients that arrived after the start of the service).

For $n = 0, \ldots, d$ define the matrices [this time they have size $(d+1) \times (d+1)$ instead of $d \times d$]

$$\bar{A}^{(n)} = \begin{pmatrix} \dfrac{\lambda_{[n+1]} + \gamma_1 \mu}{\mu} & 1 & 0 & 0 & \ldots & 0 & 0 \\ \dfrac{\lambda_{[n+2]} + \gamma_2 \mu}{\mu} & 0 & 1 & 0 & \ldots & 0 & 0 \\ \dfrac{\lambda_{[n+3]} + \gamma_3 \mu}{\mu} & 0 & 0 & 1 & \ldots & 0 & 0 \\ \vdots & \vdots & \vdots & \vdots & \ddots & \vdots & \vdots \\ \dfrac{\lambda_{[n+d-1]} + \gamma_{d-1} \mu}{\mu} & 0 & 0 & 0 & \ldots & 1 & 0 \\ \dfrac{\lambda_{[n+d]} + \gamma_d \mu}{\mu} & 0 & 0 & 0 & \ldots & 0 & 1 \\ \dfrac{\lambda_n + \gamma_0 \mu}{\mu} & 0 & 0 & 0 & \ldots & 0 & 0 \end{pmatrix},$$

where $(\mu, \gamma) \sim \nu_n$, and let $\bar{A} = \bar{A}^{(d)} \ldots \bar{A}^{(0)}$, where $\bar{A}^{(d)}, \ldots, \bar{A}^{(0)}$ are independent. Again, the method of this paper applies to this model as well, and it can be shown that Theorem 1.1 remains valid.

**2. Proof of Theorem 1.1.** The remaining part of this paper is organized as follows. First, in Section 2.1, we prove two technical lemmas and recall a result from [2] that will be used for establishing the existence of moments of $\tau$. In Section 2.2 we define the model in a formal way and introduce some more notation. In Section 2.3 we define a fluid model related to our queueing system and prove various technical facts related to it. Then, in Section 2.4 we prove part (i), in Section 2.5 we prove parts (ii) and (iii), and in Section 2.6 we prove part (iv) of Theorem 1.1.

2.1. *Some preliminary facts.* In this section we establish some technical facts needed in the course of the proof of Theorem 1.1.

Let $\Theta_1, \Theta_2, \Theta_3, \ldots$ be an i.i.d. sequence of nonnegative random matrices, and define for any $s > 0$

(2.1) $$\mathsf{k}^{\Theta}(s) = \lim_{n \to \infty} (\mathbb{E} \|\Theta_n \ldots \Theta_1\|^s)^{1/n}.$$

The following result will be important in the course of the proof of Theorem 1.1. For the scalar case, related results can be found in [4, 10], and the recent paper [8] deals with the matrix case. In fact, the above mentioned papers contain some deep results on the distribution of the random matrix (or variable) $\Lambda$ below; since the proof of Lemma 2.1 is short and elementary, we nevertheless include it to keep the paper self-contained.

LEMMA 2.1. *Abbreviate* $\Lambda = \Theta_1 + \Theta_2\Theta_1 + \Theta_3\Theta_2\Theta_1 + \cdots$. *Then:*

(i) *if* $\mathsf{k}^\Theta(s) > 1$, *then* $\mathbb{E}\|\Lambda\|^s = \infty$;
(ii) *if* $\mathsf{k}^\Theta(s) < 1$, *then* $\mathbb{E}\|\Lambda\|^s < \infty$.

PROOF. The proof of part (i) is straightforward: since the matrices are a.s. nonnegative, we have for any $n$

$$\mathbb{E}\|\Lambda\|^s \geq \mathbb{E}\|\Theta_n \ldots \Theta_1\|^s.$$

If $\mathsf{k}^\Theta(s) > 1$, then $\mathbb{E}\|\Theta_n \ldots \Theta_1\|^s \to \infty$ as $n \to \infty$ by (2.1), so (i) is proved.

Now, we turn to part (ii). Since for any $s \leq 1$

$$\|\Lambda\|^s \leq \left(\sum_{n=1}^\infty \|\Theta_n \ldots \Theta_1\|\right)^s$$
$$\leq \sum_{n=1}^\infty \|\Theta_n \ldots \Theta_1\|^s,$$

in this case the proof immediately follows from the definition of $\mathsf{k}^\Theta(s)$ as well. So, from now on we suppose that $s > 1$. We use here some technique similar to the proof of Theorem 3.1 of [11], adapted to the matrix case. Let us recall the following elementary consequence of the Jensen inequality: for any real numbers $u, v, G, H > 0$ and $s > 1$, we have

$$\frac{uG + vH}{u + v} \leq \left(\frac{uG^s + vH^s}{u + v}\right)^{1/s},$$

or, equivalently,

(2.2) $$(uG + vH)^s \leq (u + v)^{s-1}(uG^s + vH^s).$$

Denote $\Lambda_n = \Theta_1 + \Theta_2\Theta_1 + \cdots + \Theta_n \ldots \Theta_1$, and let $a_n = \mathbb{E}\|\Lambda_n\|^s$. Fix any $\alpha \in ((\mathsf{k}^\Theta(s))^{1/s}, 1)$, and apply (2.2) with $u = 1$, $v = \alpha^{n+1}$, $G = \|\Lambda_n\|$, $H = \alpha^{-(n+1)}\|\Theta_{n+1} \ldots \Theta_1\|$ to obtain that

(2.3)
$$\mathbb{E}\|\Lambda_{n+1}\|^s \leq \mathbb{E}(\|\Lambda_n\| + \alpha^{n+1} \times \alpha^{-(n+1)}\|\Theta_{n+1} \ldots \Theta_1\|)^s$$
$$\leq (1 + \alpha^{n+1})^{s-1}(\mathbb{E}\|\Lambda_n\|^s + \alpha^{n+1}\alpha^{-s(n+1)}\mathbb{E}\|\Theta_{n+1} \ldots \Theta_1\|^s).$$

From the definition of $\mathsf{k}^\Theta(s)$ it follows that

$$\lim_{n \to \infty} \alpha^{-s(n+1)}\mathbb{E}\|\Theta_{n+1} \ldots \Theta_1\|^s = 0$$

for any $\alpha$ such that $(\mathsf{k}^\Theta(s))^{1/s} < \alpha < 1$. Since, trivially, $a_n = \mathbb{E}\|\Lambda_n\|^s \geq \mathbb{E}\|\Theta_1\|^s > 0$ for any $n$, we can conclude that there exists $n_0$ such that for all $n \geq n_0$ it holds that

$$\alpha^{-s(n+1)}\mathbb{E}\|\Theta_{n+1} \ldots \Theta_1\|^s \leq a_n.$$

So, (2.3) implies that, for all $n \geq n_0$,

$$a_{n+1} \leq (1 + \alpha^{n+1})^{s-1}(a_n + \alpha^{n+1}a_n) = (1 + \alpha^{n+1})^s a_n.$$

Thus, for any $n$,

(2.4) $$\mathbb{E}\|\Lambda_n\|^s \leq M,$$

where

$$M = a_{n_0} \prod_{m=n_0}^{\infty} (1 + \alpha^{m+1})^s < \infty.$$

The proof of part (ii) of Lemma 2.1 for the case $s > 1$ follows now from (2.4) and the monotone convergence theorem. □

The next result is very similar to Lemma 2.1, so we only outline its proof.

LEMMA 2.2. *Let $X_0, X_1, X_2, \ldots$ be a sequence of a.s. positive identically distributed random variables. Suppose that $\mathbb{E} X_0^s < \infty$ for a given $s > 0$. Then, for any $\alpha > 1$,*

(2.5) $$\mathbb{E}\left(X_0 + \frac{X_1}{\alpha} + \frac{X_2}{\alpha^2} + \frac{X_3}{\alpha^3} + \cdots\right)^s < \infty.$$

PROOF. As before, in the case $s \leq 1$ the proof is straightforward. When $s > 1$, the proof can be done quite analogously to the proof of Lemma 2.1(ii) [use (2.2) with $u = 1$, $v = \alpha^{-(n+1)}$, $G = X_0 + \alpha^{-1}X_1 + \cdots + \alpha^{-n}X_n$, $H = X_{n+1}$]. □

The next result is Theorem 1 of [2]:

PROPOSITION 2.1. *Let $\alpha$ be some positive real number. Suppose that we are given a $\{\mathcal{F}_n\}$-adapted stochastic process $Z_n$, $n \geq 0$, taking values in an unbounded subset of $\mathbb{R}_+$. Denote by $\tau_\alpha$ the moment when the process $Z_n$ enters the set $(0, \alpha)$. Assume that there exist $H > 0$, $p_0 \geq 1$ such that, for any $n$, $Z_n^{p_0}$ is integrable and*

(2.6) $$\mathbb{E}(Z_{n+1}^{p_0} - Z_n^{p_0} \mid \mathcal{F}_n) \leq H Z_n^{p_0 - 1}$$

*on $\{\tau_\alpha > n\}$. Then there exists a positive constant $C = C(H, p_0)$ such that, for all $x \geq 0$ whenever $Z_0 = x$ with probability 1,*

(2.7) $$\mathbb{E}\tau_\alpha^{p_0} \leq C x^{p_0}.$$

2.2. *Notation and the formal definition of the process.* According to the description of the process given in Section 1.1, the current state of the system can be described by a pair $(R, S)$, where $R = (\mu, \gamma)$ is the current service regime, and $S = (N, \xi_0, \ldots, \xi_d)$ is the current configuration, meaning that the server is at station $N$, and $\xi_i$ is the number of customers at station $[N+i]$, $i = 0, \ldots, d$. To recover the values of $N$ and $(\xi_0, \ldots, \xi_d)$ from $S$, we use the notation $\mathcal{N}(S) = N$ and $\mathcal{K}(S) = (\xi_0, \ldots, \xi_d)$.

Denote by $\varphi_m$ the operation of the cyclic shift by $m$ positions, that is, $\varphi_m(x_0, \ldots, x_d) = (x_m, \ldots, x_d, x_0, \ldots, x_{m-1})$, and let $e_k$ be the $k$th unit vector, that is, $e_k = (0, \ldots, 0, 1, 0, \ldots, 0)$ with 1 on $k$th position, $k = 0, \ldots, d$. The transition rates are then described as follows [the current regime is $(\mu, \gamma)$]:

- $(N, \xi) \to (N, \xi + e_i)$ with rate $\lambda_{[N+i]}$, $i = 0, \ldots, d$;
- if $\xi_0 \geq 2$, then
  - $(N, \xi) \to (N, \xi - e_0 + e_i)$ with rate $\mu \gamma_i$;
  - $(N, \xi) \to (N, \xi - e_0)$ with rate $\mu(1 - \gamma_1 - \cdots - \gamma_d)$;
- suppose that $\xi = (1, 0, \ldots, 0, \xi_k, \ldots, \xi_d)$ with $\xi_k \geq 1$. Then
  - $(N, \xi) \to ([N+k], \varphi_k(\xi - e_0 + e_i))$ with rate $\mu \gamma_i$ for $i \geq k$, and the new regime is chosen independently according to $\nu_{[N+k]}$;
  - $(N, \xi) \to ([N+i], \varphi_i(\xi - e_0 + e_i))$ with rate $\mu \gamma_i$ for $i < k$, and the new regime is chosen independently according to $\nu_{[N+i]}$;
  - $(N, \xi) \to ([N+k], \varphi_k(\xi - e_0))$ with rate $\mu(1 - \gamma_1 - \cdots - \gamma_d)$, and the new regime is chosen independently according to $\nu_{[N+k]}$;
- $(N, 1, 0, \ldots, 0) \to (\mathfrak{S}, 0, \ldots, 0)$ with rate $\mu(1 - \gamma_1 - \cdots - \gamma_d)$;
- $(\mathfrak{S}, 0, \ldots, 0) \to (n, 1, 0, \ldots, 0)$ with rate $\lambda_n$, $n = 0, \ldots, d$.

We refer to the state of the process at time $t$ as $(R(t), S(t))$. Recall that we suppose that the process starts from some nonzero configuration [zero configuration is $(\mathfrak{S}, 0, \ldots, 0)$], and $\tau$ is the time of reaching the zero configuration, that is,

$$\tau = \inf\{t > 0 : \mathcal{N}(S(t)) = \mathfrak{S}\}.$$

For definiteness, we suppose that the trajectories are right-continuous and $\mathcal{N}(S(0)) = 0$. Now, we need to define the sequence of times when the server jumps. Put $t'_0 := 0$, and

$$t'_{n+1} = \inf\{t \geq t'_n : \mathcal{N}(S(t)) \neq \mathcal{N}(S(t'_n))\}$$

for all $n \geq 0$. Define also, for $n$ such that $t'_n < \tau$,

$$k^{(n)} = [\mathcal{N}(S(t'_n)) - \mathcal{N}(S(t'_{n-1}))],$$

$n \geq 1$, and $k^{(0)} := 0$. Then, define $\sigma_0 = t'_0 = 0$, and $\sigma_i = t'_n$ for all $i$ such that $k^{(0)} + \cdots + k^{(n-1)} < i \leq k^{(0)} + \cdots + k^{(n)}$. This means that if the server jumps, say, from $i$ to $[i+k]$ (i.e., it empties the queue in station $i$ and there are no customers in $[i+1], \ldots, [i+k-1]$), then we visualize this jump as an instantaneous sequence of $k$ jumps. This approach has the advantage that now we can write

$$\mathcal{N}(S(\sigma_n)) = [n] \tag{2.8}$$

for all $n$ such that $\sigma_n < \tau$.

Since we choose the regimes independently of the past evolution of the process, one also can use the following construction. Let $\mathcal{R} = (R_0, R_1, R_2, \ldots)$ be a sequence of independent $\mathcal{M}$-valued random variables, $R_i \sim \nu_{[i]}$. We can first fix a realization of this sequence, and then start the process, and the regimes will be chosen in the following way. Suppose that the initial configuration is $S(0) \neq 0$; then take $R(0) = R_{\mathcal{N}(S(0))}$ [in most cases we suppose that $\mathcal{N}(S(0)) = 0$]. Inductively, if $R(t) = R_i$ for $t \in [t'_n, t'_{n+1})$, then $R(t'_{n+1}) = R_{i+k^{(n+1)}}$. Note that, analogously to (2.8) when the server uses the regime $R_i$, it is in the station $[i]$ (at least before $\tau$; it is quite easy to modify this construction so that the last statement would be valid for all $t$, but we really do not need that for the purposes of Theorem 1.1).

In what follows, we use symbols $\mathsf{P}_\mathcal{R}$ and $\mathsf{E}_\mathcal{R}$ to denote probability and expectation given the sequence of regimes $\mathcal{R}$.

2.3. *Random fluid model and its relationship with the original process.* Now, we define what we call the *random fluid model.* Here we suppose that $\lambda$'s and $\mu$ are flow rates, that is, the liquid comes to node $n$ with speed $\lambda_n$, and when the server is in station $k$ with service rate and feedback parameters $(\mu, \gamma)$, it sends the liquid from $k$ out of the system with speed $\mu(1 - \gamma_1 - \cdots - \gamma_d)$, and to each site $[k+j]$ with speed $\mu\gamma_j$, $j = 1, \ldots, d$. In other words, the randomness present in this model only relates to the fact that, when arriving to a node, the server chooses the serving and transition parameters (i.e., the speeds of the corresponding flows) at random.

We will refer to $(R^f(t), S^f(t))$ as the state of the random fluid model at time $t$ [where, as before, $S^f(t)$ is the configuration and $R^f(t)$ is the regime]. If for a state $(R^f, S^f)$ the regime is $R^f = (\mu, \gamma)$ and the configuration is $S^f = (N, x_0, \ldots, x_d)$, then we define by

$$\mathfrak{T}(R^f, S^f) := \frac{x_0}{\mu - \lambda_N}$$

the time that the server needs to pump all the liquid out of the (current) station $N$.

To define the fluid model starting from $S^f(0) = (N, x_0^0, \ldots, x_d^0)$ (as mentioned before, normally we will take $N = 0$), where $x_0^0 > 0$, and with the initial regime $R^f(0) = (\mu, \gamma)$, put

$$S^f(t) = (N, x_0^0 - (\mu - \lambda_N)t,$$
$$x_1^0 + (\mu\gamma_1 + \lambda_{[N+1]})t, \ldots, x_d^0 + (\mu\gamma_d + \lambda_{[N+d]})t),$$

and $R^f(t) = R^f(0)$ for $t < \sigma_{N+1}^f := \mathfrak{T}(R^f(0), S^f(0)) = x_0^0/(\mu - \lambda_N)$. For $t = \sigma_{N+1}^f$, choose $R^f(\sigma_{N+1}^f) \sim \nu_{[N+1]}$ independently, and put

$$S^f(\sigma_{N+1}^f) = ([N+1], x_1^0 + (\mu\gamma_1 + \lambda_{[N+1]})\sigma_{N+1}^f, \ldots,$$
$$x_d^0 + (\mu\gamma_d + \lambda_{[N+d]})\sigma_{N+1}^f, 0).$$

Inductively, suppose that for some $i > N$ the server switching time $\sigma_i^f$ is defined, and the current configuration is $S^f(\sigma_i^f) = ([i], x_0^i, \ldots, x_{d-1}^i, 0)$ (note that the last coordinate is 0 because the server has just emptied the station and jumped, and that $\mathcal{N}(S^f(\sigma_i^f))$ always equals $[i]$ since, for the fluid model, the server can never jump over a station). Then, choose $R^f(\sigma_i^f) \sim \nu_{[i]}$ independently, and supposing that $R^f(\sigma_i^f) = (\mu', \gamma')$, put $\sigma_{i+1}^f = \sigma_i^f + \mathfrak{T}(R^f(\sigma_i^f), S^f(\sigma_i^f))$, for $t \in [\sigma_i^f, \sigma_{i+1}^f)$ define $R^f(t) = R^f(\sigma_i^f)$:

$$S^f(t) = ([i], x_0^i - (\mu' - \lambda_{[i]})t,$$
$$x_1^i + (\mu'\gamma_1' + \lambda_{[i+1]})t, \ldots, x_d^i + (\mu'\gamma_d' + \lambda_{[i+1]})t),$$

and [abbreviating for the moment $\mathfrak{T} := \mathfrak{T}(R^f(\sigma_i^f), S^f(\sigma_i^f))$]

$$S^f(\sigma_{i+1}^f) = ([i+1], x_1^i + (\mu'\gamma_1' + \lambda_{[i+1]})\mathfrak{T}, \ldots, x_d^i + (\mu'\gamma_d' + \lambda_{[i+1]})\mathfrak{T}, 0).$$

Let $S^f(\sigma_{i+1}^f) := ([i+1], x_0^{i+1}, \ldots, x_{d-1}^{i+1}, 0)$. At this point it is important to note that, denoting by $x^i$ (resp., $x^{i+1}$) the column vector with coordinates $x_j^i$ (resp., $x_j^{i+1}$), $j = 0, \ldots, d-1$, we have

$$(2.9) \qquad x^{i+1} = \tilde{A}(i)x^i,$$

where $\tilde{A}(i)$ is the matrix defined in (1.1), with $(\mu', \gamma')$ instead of $(\mu, \gamma)$ [i.e., $\tilde{A}(i)$ depends on the current regime]. Letting $k = [i]$, $(\mu', \gamma')$ was sampled from $\nu_k$ and $\tilde{A}(i)$ is a random matrix distributed like $A^{(k)}$ of (1.1).

An important technical fact is that, for the fluid model, at the moments $\sigma_i^f$ of server's transition, the relative difference between the amounts of liquid in the stations (not counting the station that the server just left) cannot become too large:

LEMMA 2.3. *Suppose that the initial configuration $S^f(0) = (N, x_0^0, \ldots, x_d^0)$ is such that $x_k^0 \neq 0$ for all $k = 0, \ldots, d-1$. Then there exists a constant $K \in (0, \infty)$, depending only on $\varepsilon_0, M_0$ (from Condition 1.1), $\lambda_k$'s, $d$, such that a.s. for all $i > d$,*

$$(2.10) \qquad \frac{x_n^i}{x_m^i} \leq K \qquad \text{for all } n, m = 0, \ldots, d-1.$$

PROOF. Denote $\lambda = \min\{\lambda_0,\ldots,\lambda_d\} > 0$, $\hat{\lambda} = \max\{\lambda_0,\ldots,\lambda_d\}$, and suppose without restriction of generality that $N = 0$. Let $h = \sum_{i=0}^{d} x_i^0$ be the total amount of liquid present in the fluid system at time 0. Note that, since we aim for a constant $K$ that does not depend on the initial configuration, to prove this lemma, it is enough to show that there is $K$ such that

$$\frac{\max\{x_0^{d+1},\ldots,x_{d-1}^{d+1}\}}{\min\{x_0^{d+1},\ldots,x_{d-1}^{d+1}\}} \leq K \qquad \text{a.s.} \tag{2.11}$$

First, we will give a lower bound on $\min\{x_0^{d+1},\ldots,x_{d-1}^{d+1}\}$. Note that, according to Condition 1.1, each station is emptied with speed at most $M_0$, a.s. Clearly, $\max\{x_0^{d+1},\ldots,x_d^{d+1}\} \geq h/(d+1)$, and suppose that this maximum is reached in station $\hat{m}$. Then, the time it takes to empty the $\hat{m}$th station $\mathfrak{T}(R^f(\sigma_{\hat{m}}^f), S^f(\sigma_{\hat{m}}^f))$ is at least $\frac{h}{(d+1)M_0}$ (since the amount of liquid in $\hat{m}$ will not decrease by the time the server gets there). Therefore, at the moment $\sigma_{\hat{m}+1}^f$ (i.e., when the server pumps all the liquid out of station $\hat{m}$), we obtain that *all* the nonempty stations have at least $\frac{\lambda h}{(d+1)M_0}$ units of liquid. Repeating this argument once again (if $\hat{m} \neq d$), we obtain that

$$\min\{x_0^{d+1},\ldots,x_{d-1}^{d+1}\} \geq \frac{\lambda^2 h}{(d+1)M_0^2}. \tag{2.12}$$

Now, we obtain an upper bound on $\max\{x_0^{d+1},\ldots,x_{d-1}^{d+1}\}$. Again, the key observation is that, by Condition 1.1, the speed of emptying the stations cannot be less than $\varepsilon_0$, a.s. Initially, each station (including the one where the server starts) has at most $h$ units of liquid. Thus, the time $\mathfrak{T}(R^f(0), S^f(0))$ to empty the 0th station is at most $h/\varepsilon_0$, and so at time $\sigma_1^f$ all the stations have at most

$$h + \frac{(\hat{\lambda}+1)h}{\varepsilon_0} = h\left(1 + \frac{(\hat{\lambda}+1)}{\varepsilon_0}\right)$$

units of liquid. Repeating this argument $d$ times more, we see that

$$\max\{x_0^{d+1},\ldots,x_{d-1}^{d+1}\} \leq h\left(1 + \frac{(\hat{\lambda}+1)}{\varepsilon_0}\right)^{d+1}. \tag{2.13}$$

Using (2.12) and (2.13), we find that (2.11) is verified with

$$K = \frac{(1 + \varepsilon_0^{-1}(\hat{\lambda}+1))^{d+1}(d+1)M_0^2}{\lambda^2},$$

and the proof of Lemma 2.3 is finished. $\square$

Define a sequence of random variables

$$D_n(S^f(0)) = \sigma_n^f = \sum_{i=0}^{n-1} \mathfrak{T}(R^f(\sigma_i^f), S^f(\sigma_i^f))$$

(where $\sigma_0^f := 0$), and

$$D_\infty(S^f(0)) = \lim_{n \to \infty} \sigma_n^f = \sum_{i=0}^{\infty} \mathfrak{T}(R^f(\sigma_i^f), S^f(\sigma_i^f)).$$

Intuitively, the random variable $D_\infty(S^f(0))$ can be interpreted as the time needed to completely empty the fluid system. Clearly, when $D_\infty(S^f(0)) < \infty$, the above definition of the fluid model works only when $t$ is less than the (random) time $D_\infty(S^f(0))$.

Recall that, by Condition 1.1, it holds that

(2.14) $$\mathfrak{T}(R^f(\sigma_i^f), S^f(\sigma_i^f)) \in [M_0^{-1} x_0^i, \varepsilon_0^{-1} x_0^i] \quad \text{a.s.}$$

Also, from Condition 1.1 and (2.9), we obtain that for some positive constants $C_1, C_2$ and for all $k$

(2.15) $$C_1 \|x^k\| \leq \|x^{k+1}\| \leq C_2 \|x^k\| \quad \text{a.s.}$$

Using (2.14) and (2.15), we obtain that for all $i \geq 0$

(2.16) $$\sum_{j=0}^{d} \mathfrak{T}(R^f(\sigma_{i+j}^f), S^f(\sigma_{i+j}^f)) \leq \varepsilon_0^{-1}(x_0^i + \cdots + x_0^{i+d})$$
$$\leq \varepsilon_0^{-1}(1 + C_2 + \cdots + C_2^d)\|x^i\|.$$

Analogously, using also Lemma 2.3, we can write for all $i \geq d$

(2.17) $$\sum_{j=0}^{d} \mathfrak{T}(R^f(\sigma_{i+j}^f), S^f(\sigma_{i+j}^f))$$
$$\geq M_0^{-1}(x_0^i + \cdots + x_0^{i+d})$$
$$\geq (M_0 K d)^{-1}(1 + (C_1/Kd)^{-1} + \cdots + (C_1/Kd)^d)\|x^i\|.$$

Denote $\hat{A}(n) = \tilde{A}(d + n(d+1)) \ldots \tilde{A}(1 + n(d+1))\tilde{A}(n(d+1))$. Clearly, the matrices $\hat{A}(n)$ are i.i.d. and $\hat{A}(n) \stackrel{\text{law}}{=} A$. Denote also for all $0 \leq n \leq \infty$

$$\hat{\Lambda}_n = \hat{A}(0) + \hat{A}(1)\hat{A}(0) + \cdots + \hat{A}(n) \ldots \hat{A}(0).$$

Note the following two properties of the $\mathcal{L}_1$-norm:

- if $B_1, B_2$ are nonnegative matrices and $x$ is a nonnegative vector, then $\|B_1 x\| + \|B_2 x\| = \|(B_1 + B_2)x\|$, and
- if $\|x\| \leq 1$ and $x_j \geq c > 0$ for all $j$, then $\|Bx\| \leq \|B\| \leq c^{-1}\|Bx\|$ for any nonnegative matrix $B$.

Taking this into account, from (2.9), (2.16) and (2.17), we obtain that there exist $U_1, U_2$ such that

$$U_1 \|\hat{\Lambda}_{n-1}\| \leq D_{n(d+1)}(S^f(0)) \leq U_2 \|\hat{\Lambda}_{n-1}\|, \tag{2.18}$$

and so

$$U_1 \|\hat{\Lambda}_\infty\| \leq D_\infty(S^f(0)) \leq U_2 \|\hat{\Lambda}_\infty\|. \tag{2.19}$$

Analogously to the stochastic model of Section 2.2, we can first choose the sequence of regimes $\mathcal{R} = (R_0, R_1, R_2, \ldots)$ (as before, it is a sequence of independent $\mathcal{M}$-valued random variables, $R_i \sim \nu_{[i]}$). Then, given a realization of this sequence, put $R^f(t) := R_i$, $t \in [\sigma_i^f, \sigma_{i+1}^f)$, $i \geq \mathcal{N}(S^f(0))$. It is important to note that, when $\mathcal{R}$ is fixed, the fluid model becomes a completely deterministic process; in particular, the values of $D_n(S^f(0))$ and $D_\infty(S^f(0))$ can be explicitly calculated. To indicate this, when $\mathcal{R}$ is given, we write $D_n(S^f(0) \mid \mathcal{R})$ and $D_\infty(S^f(0) \mid \mathcal{R})$ for those quantities. Also, note that all the variables in $\mathcal{R}$ [beginning with the $\mathcal{N}(S^f(0))$th one] will be used, which is not necessarily the case for the stochastic model.

To prove Theorem 1.1, an important idea is to relate the behavior of the stochastic system to that of the fluid model when they share the same regimes. In other words, a coupling of these two systems can be performed by using the same sequence of regimes $\mathcal{R}$. Consider any nonempty initial configuration $S(0) = S^f(0) = (n, y_0, \ldots, y_d)$ and take any $R(0) = R^f(0) = (\mu, \gamma) \in \operatorname{supp} \nu_n$. Let $\sigma^s = \inf\{t \geq 0 : \mathcal{N}(S(t)) \neq n\}$ be the time one needs to wait until the server's jump for the stochastic model, and let $\mathfrak{T}(R^f(0), S^f(0)) = \frac{y_0}{\mu - \lambda_n}$ be the corresponding time for the fluid model. Then, we have the following result:

LEMMA 2.4. (i) *For any $\delta > 0$ there exists a constant $V_1 = V_1(\delta) > 0$ (depending also on $\varepsilon_0, M_0$ from Condition 1.1, $\lambda$s, d) such that*

$$\mathsf{P}_\mathcal{R}[|\sigma^s - \mathfrak{T}(R^f(0), S^f(0))| \leq \delta y_0] \geq 1 - e^{-V_1 y_0}. \tag{2.20}$$

(ii) *Suppose that $S(\sigma^s) = ([n+1], z_0, \ldots, z_{d-1}, 0)$ (here we assume that in the initial configuration $y_1 > 0$, otherwise it is not certain that the next destination of the server will be $[n+1]$) and $S^f(\mathfrak{T}(R^f(0), S^f(0))) = ([n+1], z_0^f, \ldots, z_{d-1}^f, 0)$. Let $\tilde{A}$ be the matrix defined in (1.1), and abbreviate $y = (y_0, \ldots, y_{d-1})^\mathsf{T}$, $z = (z_0, \ldots, z_{d-1})^\mathsf{T}$, $z^f = (z_0^f, \ldots, z_{d-1}^f)^\mathsf{T} = \tilde{A}y$. Then for any $\delta > 0$ there exist $V_2 = V_2(\delta) > 0$ (depending also on $\varepsilon_0, M_0$ from Condition 1.1, $\lambda_k$'s, d) such that*

$$\mathsf{P}_\mathcal{R}[|z - z^f| \leq \delta y_0] \geq 1 - e^{-V_2 y_0}. \tag{2.21}$$

PROOF. The proof of this fact is completely elementary, so we give only a sketch.

First, we deal with (2.20). The evolution of the queue in the current station $n$ can be represented as a birth-and-death process $X(t)$ on $\mathbb{Z}$ starting from $y_0 > 0$, with (constant) birth rate $\lambda_n$ and death rate $\mu$. Now, the moment $\sigma^s$ is the time to reach 0 for this process. For the sake of simplicity, suppose that $\mu + \lambda_n = 1$, that is, the total jump rate of the process $X(t)$ is 1 (this can always be done by rescaling the time, and the scaling factor will be bounded from both sides by Condition 1.1). Note that

$$\mathbb{P}\left[\sigma^s > \frac{y_0}{\mu - \lambda_n} + \delta y_0\right] \leq \mathbb{P}\left[X\left(\frac{y_0}{\mu - \lambda_n} + \delta y_0\right) > 0\right].$$

The number of jumps of $X(t)$ on the time interval $[0, \frac{y_0}{\mu - \lambda_n} + \delta y_0]$ will be at least $\frac{y_0}{\mu - \lambda_n} + \frac{\delta y_0}{2}$ with probability at least $1 - e^{-C_1 y_0}$ for some $C_1 > 0$ (and it is not difficult to see that the constant $C_1$ can be chosen uniformly by Condition 1.1). For the corresponding discrete-time random walk $\hat{X}(n)$, it is a straightforward computation to bound $\mathbb{P}[\hat{X}(\frac{y_0}{\mu - \lambda_n} + \frac{\delta y_0}{2}) > 0]$ from above.

To bound $\mathbb{P}[\sigma^s < \frac{y_0}{\mu - \lambda_n} - \delta y_0]$ from above, we note that, by time $\frac{y_0}{\mu - \lambda_n} - \delta y_0$, the continuous-time process $X(t)$ will perform at most $\frac{y_0}{\mu - \lambda_n} - \frac{\delta y_0}{2}$ jumps with probability at least $1 - e^{-C_2 y_0}$ for some $C_2 > 0$. So, we write

$$\mathbb{P}\left[\sigma^s < \frac{y_0}{\mu - \lambda_n} - \delta y_0\right] \leq \sum_{m \leq (y_0/(\mu - \lambda_n)) - \delta y_0/2} \mathbb{P}[\hat{X}(m) \leq 0],$$

and then again use large deviations estimates to bound the terms in the right-hand side of the above display.

As for (2.21), it follows easily from (2.20) if we take into account the following two observations:

- by (2.20), we have good control on $\sigma^s$ (the time it takes to empty the queue in $n$);
- knowing this time, for all $i = 1, \ldots, d$ it is straightforward to estimate how many customers will come to station $[n+i]$ during that time (given $\sigma^s = t$, the number of customers that came from outside is simply a Poisson random variable with mean $\lambda_{[n+i]}t$); moreover, we can easily write estimates on the total number of customers that came out of station $n$, and then observe that the proportion (roughly) $\gamma_i$ of them go to $[n+i]$. □

2.4. *Transience.* Here we prove part (i) of Theorem 1.1, that is, we show that if $\mathsf{k}'(0) > 0$, then $\mathbb{P}[\tau = \infty] > 0$. The corresponding result in [11] was established using the Lyapunov function technique, but we provide a different argument here.

The quantity $\mathsf{k}'(0) > 0$ is in fact the top Lyapunov exponent of $A$, so for almost every $\mathcal{R}$ there exist $\alpha > 1$ (not depending on $\mathcal{R}$) and a positive number $C_1 = C_1(\mathcal{R})$ such that

$$\|A_n \ldots A_1\| > C_1 \alpha^n \tag{2.22}$$

for all $n$. Note that, since $A \stackrel{\text{law}}{=} \tilde{A}(i(d+1)+d)\ldots\tilde{A}(i(d+1))$ for all $i \geq 0$, from (2.22) we obtain that for some $C_2$

$$\|\tilde{A}(m-1)\ldots\tilde{A}(d+1)\| \geq C_2 \alpha^{m/(d+1)} \tag{2.23}$$

for all $m > d+1$.

Now, the idea is to compare the stochastic model and the fluid model which start from the same initial configuration and share the same sequence of regimes.

Consider the stochastic and the fluid model with the same initial configuration:

$$S(0) = S^f(0) = (0, \xi_0^0, \ldots, \xi_{d-1}^0, 0) = (0, x_0^0, \ldots, x_{d-1}^0, 0),$$

and, as before, $(\xi_0^i, \ldots, \xi_d^i) = \mathcal{K}(S(\sigma_i))$, $(x_0^i, \ldots, x_d^i) = \mathcal{K}(S^f(\sigma_i))$. Since we have $\xi_d^i = x_d^i = 0$ for all $i$, let us denote $\xi^i = (\xi_0^i, \ldots, \xi_{d-1}^i)^\mathsf{T}$ and $x^i = (x_0^i, \ldots, x_{d-1}^i)^\mathsf{T}$. Here $\mathsf{T}$ means transposed, so $\xi^i$ and $x^i$ are column vectors.

Consider the sequence of events

$$B_m^\delta = \{\xi^m \geq (1-\delta)^m x^m\},$$

$m \geq 1$. From (2.9) and (2.23) we obtain that

$$\|x^m\| = \|\tilde{A}(m-1)\tilde{A}(m-2)\ldots\tilde{A}(d+1)x^{d+1}\|$$
$$\geq C_3 \|\tilde{A}(m-1)\tilde{A}(m-2)\ldots\tilde{A}(d+1)\|,$$

so, by Lemma 2.3, we have that $x_0^m \geq C_4 \alpha^{m/(d+1)}$ for some $C_4 > 0$. With this fact, we recall Lemma 2.4(ii), and write for $n > d+1$

$$\begin{aligned}\mathsf{P}_\mathcal{R}&[B_n^\delta \mid B_{n-1}^\delta, \ldots, B_1^\delta] \\ &\geq \mathbb{P}[\xi^n \geq (1-\delta)^{n-1}\tilde{A}(n-1)x^{n-1} \mid \xi^{n-1} = (1-\delta)^{n-1}x^{n-1}] \\ &\geq 1 - \exp\{-V_2(\delta)(1-\delta)^{n-1}x_0^{n-1}\} \\ &\geq 1 - \exp\{-V_2(\delta)(1-\delta)^{n-1}C_4\alpha^{(n-1)/(d+1)}\}.\end{aligned} \tag{2.24}$$

Choose $\delta > 0$ in such a way that $(1-\delta)\alpha^{1/(d+1)} > 1$. It is clear that $\mathsf{P}_\mathcal{R}[B_1^\delta] > 0$ and $\mathsf{P}_\mathcal{R}[B_k^\delta \mid B_{k-1}^\delta, \ldots, B_1^\delta] > 0$ for all $k \leq d+1$, so, by (2.24), it follows that $\mathsf{P}_\mathcal{R}[\bigcap_{n=1}^\infty B_n^\delta] > 0$. Since on $\{\bigcap_{n=1}^\infty B_n^\delta\}$ we know that $\xi_0^n \geq (1-\delta)^n C_4 \times \alpha^{n/(d+1)} \to \infty$, it holds that $\tau = \infty$ with positive probability for almost every $\mathcal{R}$. This proves the transience result.

2.5. *Recurrence and existence of moments.* First, note that if $\mathsf{k}'(0) < 0$, then $s_0 > 0$, and that if $\mathbb{E}\tau^s < \infty$ for at least one $s > 0$, then $\mathbb{P}[\tau = \infty] = 0$, so part (ii) of Theorem 1.1 easily follows from part (iii). Thus, we concentrate on part (iii).

Fix $s$ such that $\mathsf{k}(s) < 1$ (or, equivalently, $s < s_0$). For the stochastic process $S(t)$ define a sequence of stopping times $\ell_n$, $n = 0, 1, 2, \ldots$, via $\ell_0 := 0$,

$$\ell_{n+1} = \inf\{t > \ell_n : S(t) \neq S(\ell_n)\}$$

for $n \geq 0$. The time $\ell_n$ is the $n$th jump time of $S(t)$, that is, a customer arrives at the system or a service is completed.

Given the sequence of regimes $\mathcal{R}$, the process $\hat{S}(n) := S(\ell_n)$ is a discrete-time Markov chain (however, in general it need not be time-homogeneous). If the current regime is $(\mu, \gamma)$, then if $\hat{S}(n) = (N, \xi)$ with $\xi_0 > 1$ (it is easy also to treat the case $\xi_0 = 1$ analogously to what was done in the beginning of Section 2.2 for continuous time), then

$$\hat{S}(n+1) = \begin{cases} (N, \xi + e_i) & \text{with probability } \mathcal{Z}^{-1}\lambda_{[N+i]}, \\ (N, \xi - e_0 + e_i) & \text{with probability } \mathcal{Z}^{-1}\mu\gamma_i, \\ (N, \xi - e_0) & \text{with probability } \mathcal{Z}^{-1}\mu(1 - \gamma_1 - \cdots - \gamma_d), \end{cases}$$

where $\mathcal{Z} = \mu + \lambda_0 + \cdots + \lambda_d$.

For a configuration $S = (N, \xi)$ and a sequence of regimes $\mathcal{R}$ define $f^{\mathcal{R}}(S)$ to be the time it takes for the fluid model to arrive at 0 starting from $S$:

$$f^{\mathcal{R}}(S) = D_\infty(S \mid \mathcal{R}).$$

We need the following:

LEMMA 2.5. *There exists $\hat{\varepsilon} > 0$ such that for all $n$ and for all possible regimes*

(2.25) $\quad \mathsf{E}_{\mathcal{R}}(f^{\mathcal{R}}(\hat{S}(n+1)) - f^{\mathcal{R}}(\hat{S}(n)) \mid \hat{S}(n) = S, \mathcal{N}(S) \neq \mathfrak{S}) < -\hat{\varepsilon}.$

PROOF. First of all, it is not difficult to obtain that $f^{\mathcal{R}}(N, x)$ is a linear function of $x$. Then, if $(\mu, \gamma)$ is the current regime, write, abbreviating as before $\mathcal{Z} = \mu + \lambda_0 + \cdots + \lambda_d$,

$$\mathsf{E}_{\mathcal{R}}(f^{\mathcal{R}}(\hat{S}(n+1)) \mid \hat{S}(n) = (N, x))$$

(2.26)
$$= \mathcal{Z}^{-1} \sum_{i=0}^{d} \lambda_{[N+i]} f^{\mathcal{R}}(N, x + e_i) + \mathcal{Z}^{-1}\mu \sum_{i=1}^{d} \gamma_i f^{\mathcal{R}}(N, x - e_0 + e_i)$$
$$+ \mathcal{Z}^{-1}\mu(1 - \gamma_1 - \cdots - \gamma_d)f^{\mathcal{R}}(N, x - e_0)$$
$$= f^{\mathcal{R}}(N, x + \hat{y}),$$

where

$$\hat{y} = \mathcal{Z}^{-1}\left(\sum_{i=0}^{d} \lambda_{[N+i]} e_i + \mu \sum_{i=1}^{d} \gamma_i (e_i - e_0) - \mu(1 - \gamma_1 - \cdots - \gamma_d) e_0\right)$$

(2.27)

$$= \mathcal{Z}^{-1}\left(\sum_{i=1}^{d} (\lambda_{[N+i]} + \mu \gamma_i) e_i - (\mu - \lambda_N) e_0\right).$$

If the fluid model is now in $(N, x)$ with regime $(\mu, \gamma)$, then in $\mathcal{Z}^{-1}$ time units it will be in $(N, x + \hat{y})$, so

$$\mathsf{E}_{\mathcal{R}}(f^{\mathcal{R}}(\hat{S}(n+1)) - f^{\mathcal{R}}(\hat{S}(n)) \mid \hat{S}(n) = S, \mathcal{N}(S) \neq \mathfrak{S}) = -\mathcal{Z}^{-1} < -\hat{\varepsilon},$$

where $\hat{\varepsilon} = (\lambda_0 + \cdots + \lambda_d + M_0)^{-1}$ (recall that $M_0$ is from Condition 1.1). $\square$

From Theorem 2.1.1 of [6] we obtain that (remember that $\tau$ corresponds to the continuous-time process, and the total jump rate is bounded from below by $\lambda_0 + \cdots + \lambda_d$)

(2.28) $$\mathsf{E}_{\mathcal{R}} \tau \leq \hat{\varepsilon}^{-1}(\lambda_0 + \cdots + \lambda_d) f^{\mathcal{R}}(S(0)),$$

and, by the Jensen inequality,

(2.29) $$\mathsf{E}_{\mathcal{R}} \tau^s \leq \hat{\varepsilon}^{-s}(\lambda_0 + \cdots + \lambda_d)^s (f^{\mathcal{R}}(S(0)))^s,$$

for any $s \leq 1$. We can use (2.19) and Lemma 2.1 to get that $\mathbb{E}\tau^s < \infty$ when $\mathsf{k}(s) < 1$ (for the case $s \leq 1$).

Now, the goal is to obtain an analog of (2.29) for $s > 1$.

The idea is to use Proposition 2.1. For that, it would be nice to have a constant $C$ such that

(2.30) $$|f^{\mathcal{R}}(\hat{S}(n+1)) - f^{\mathcal{R}}(\hat{S}(n))| \leq C, \quad \text{a.s.}$$

For a fixed $\mathcal{R}$ such $C$ exists [take $C = f^{\mathcal{R}}(0, e_0 + \cdots + e_d)$], but it is not uniform, which prevents us from getting an analog of (2.29) (with a constant that does not depend on $\mathcal{R}$!) directly from Proposition 2.1.

To get around this difficulty, let us proceed in the following way. We consider another function $f_\varepsilon^{\mathcal{R}}(\cdot)$, which is the time of reaching the origin for a slightly modified fluid model that dominates the original one. For this new function, when the maximal one-step increment is large [cf. (2.30)], it happens also that the drift in the negative direction is large as well, and that permits us to apply Proposition 2.1.

Fix $\varepsilon > 0$ and consider a modification of the present model, where the new arrival rates are $\lambda_i' = \lambda_i + \varepsilon$, $i = 0, \ldots, d$. As in Section 2.3, define $D_\infty^{(\varepsilon)}(S \mid \mathcal{R})$ to be the total time needed to completely empty the fluid system with $\lambda_i'$'s in the place of $\lambda_i$'s, starting with $S$ and when the sequence of regimes $\mathcal{R}$ is fixed. And, denote $f_\varepsilon^{\mathcal{R}}(S) := D_\infty^{(\varepsilon)}(S \mid \mathcal{R})$.

Analogously to (2.26), we can write (here $\mathtt{E}_\mathcal{R}$ is the same as before, i.e., it is related to the process with $\lambda_i$'s, and *not* $\lambda'_i$'s)

$$\mathtt{E}_\mathcal{R}(f_\varepsilon^\mathcal{R}(\hat{S}(n+1)) \mid \hat{S}(n) = (N,x)) = f_\varepsilon^\mathcal{R}(N, x+\hat{y})$$
$$= f_\varepsilon^\mathcal{R}(N, x+\hat{y}_1) - \varepsilon F_N^{(\varepsilon)},$$

where $\hat{y}$ is from (2.27),

$$\hat{y}_1 = \mathcal{Z}^{-1}\left(\sum_{i=1}^d (\lambda'_{[N+i]} + \mu\gamma_i)e_i - (\mu - \lambda'_N)e_0\right)$$
$$= \hat{y} + \varepsilon \mathcal{Z}^{-1}(e_0 + \cdots + e_d),$$

and $F_N^{(\varepsilon)} := f_\varepsilon^\mathcal{R}(N, e_0 + \cdots + e_d)$. So,

$$\mathtt{E}_\mathcal{R}(f_\varepsilon^\mathcal{R}(\hat{S}(n+1)) \mid \hat{S}(n) = (N,x)) = -\mathcal{Z}^{-1} - \varepsilon \mathcal{Z}^{-1} F_N^{(\varepsilon)}$$
(2.31)
$$< -\varepsilon \mathcal{Z}^{-1} F^{(\varepsilon)},$$

where $F^{(\varepsilon)} = \min_{N=0,\ldots,d} F_N^{(\varepsilon)}$. Also, from Condition 1.1 and Lemma 2.3 it is not difficult to show that [compare with (2.30)] there exists $C_1$ such that

$$|f_\varepsilon^\mathcal{R}(\hat{S}(n+1)) - f_\varepsilon^\mathcal{R}(\hat{S}(n))| \leq C_1 F^{(\varepsilon)}, \qquad \text{a.s.}$$

Again, using Condition 1.1 and Lemma 2.3, it is straightforward to obtain that there exists $C_2$ such that, on $\hat{S}(n) = (N,x)$

$$(2.32) \qquad \frac{f_\varepsilon^\mathcal{R}(\hat{S}(n+1)) - f_\varepsilon^\mathcal{R}(N,x)}{f_\varepsilon^\mathcal{R}(N,x)} \leq \frac{C_2}{\|x\|}.$$

Abbreviate $Z_n := f_\varepsilon^\mathcal{R}(\hat{S}(n))$. Using some elementary calculus, we see that for any $C_4 > 0$ there exists $C_3 > 0$ (depending on $s$) such that

$$(1+y)^s \leq 1 + sy + C_3 y^2 \qquad \text{on } \{y \in \mathbb{R} : |y| < C_4\}.$$

So, by (2.31) and (2.32), if $x$ is such that $\frac{C_2}{\|x\|} < C_4$,

$$\mathtt{E}_\mathcal{R}(Z_{n+1}^s - Z_n^s \mid \hat{S}(n) = (N,x))$$
$$= (f_\varepsilon^\mathcal{R}(N,x))^{s-1} \mathtt{E}_\mathcal{R}\left(\left(1 + \frac{Z_{n+1} - Z_n}{Z_n}\right)^s - 1 \,\Big|\, \hat{S}(n) = (N,x)\right)$$
$$\leq (f_\varepsilon^\mathcal{R}(N,x))^{s-1}[s\mathtt{E}_\mathcal{R}(Z_{n+1} - Z_n \mid \hat{S}(n) = (N,x))$$
$$\qquad + C_3(f_\varepsilon^\mathcal{R}(N,x))^{-1} \mathtt{E}_\mathcal{R}((Z_{n+1} - Z_n)^2 \mid \hat{S}(n) = (N,x))]$$
$$\leq (f_\varepsilon^\mathcal{R}(N,x))^{s-1}\left(-s\mathcal{Z}^{-1}\varepsilon F^{(\varepsilon)} + C_3 C_1^2 \frac{(F^{(\varepsilon)})^2}{\|x\| f_\varepsilon^\mathcal{R}(N, x/\|x\|)}\right).$$

By Condition 1.1 and Lemma 2.3, there exist $C_5, C_6 > 0$ such that $F^{(\varepsilon)} > C_5$ and $\frac{F^{(\varepsilon)}}{f_\varepsilon^{\mathcal{R}}(N,e)} < C_6$ for any $e$ such that $\|e\| = 1$. So, we obtain that there exists $a_1$ such that

$$\mathsf{E}_{\mathcal{R}}(Z_{n+1}^s - Z_n^s \mid \hat{S}(n) = (N,x)) \leq (f_\varepsilon^{\mathcal{R}}(N,x))^{s-1} F^{(\varepsilon)}\left(-s\mathcal{Z}^{-1}\varepsilon + \frac{C_3 C_1^2 C_6}{\|x\|}\right)$$
$$\leq -C_7(f_\varepsilon^{\mathcal{R}}(N,x))^{s-1}$$

for $x$ such that $\|x\| \geq a_1$. From Proposition 2.1, we obtain that there exists $C_8 > 0$ (depending also on $s$) such that

$$\mathsf{E}_{\mathcal{R}}\tau^s \leq C_8(f_\varepsilon^{\mathcal{R}}(S(0)))^s, \tag{2.33}$$

for $s > 1$.

Now, for any $\delta > 0$ we can find $\varepsilon > 0$ such that

$$f_\varepsilon^{\mathcal{R}}(S(0)) = D_\infty^{(\varepsilon)}(S(0) \mid \mathcal{R}) \leq \|\hat{\Lambda}_\infty^{(\delta)}\| \|S(0)\|, \tag{2.34}$$

where

$$\hat{\Lambda}_\infty^{(\delta)} = (1+\delta)\hat{A}(0) + (1+\delta)^2 \hat{A}(1)\hat{A}(0) + (1+\delta)^3 \hat{A}(2)\hat{A}(1)\hat{A}(0) + \cdots.$$

Choose $\delta$ in such a way that $(1+\delta)\mathsf{k}(s) < 1$, and use (2.34) and Lemma 2.1 to obtain that

$$\mathbb{E}(\tau(a_1))^s < \infty,$$

where $\tau(a_1)$ is the moment of hitting the set $\{x : \|x\| < a_1\}$ for the discrete-time process. Then, by a standard argument (see, e.g., Theorem A.1 from [13]), we get that, for the discrete-time process, the $s$th moment of the hitting time of $\mathfrak{S}$ is finite, and so $\mathbb{E}\tau^s < \infty$. $\square$

2.6. *Nonexistence of moments.* In this section we prove part (iv) of Theorem 1.1. We have to prove that if $s > s_0$ [equivalently, $\mathsf{k}(s) > 1$], then $\mathbb{E}\tau^s = \infty$.

The strategy of the proof is the following:

1. We consider the fluid model and we note first that, by Lemma 2.1, the expectation of the $s$th moment of the random variable $D_\infty$ is infinite.
2. We consider a smaller random time, which is the time of reaching some "growing" set. Then, with the help of Lemma 2.2, we prove that the $s$th moment of that random time is still infinite.
3. We look now at the stochastic model when the trajectory of the fluid model is fixed (i.e., with a fixed sequence of regimes $\mathcal{R}$). Then, we prove that with a constant probability the time of reaching that finite set for the stochastic process will be greater than the corresponding time for the fluid model.

4. This fact permits us to obtain that, given $\mathcal{R}$, the $s$th moment of the time of reaching 0 for the stochastic process is (roughly) proportional to the corresponding time for the fluid model (which is deterministic, since $\mathcal{R}$ is fixed).
5. Then, we integrate and, using Lemma 2.1, conclude that $\mathbb{E}\tau^s = \infty$.

Let $n_t = \max\{i : \sigma_i^f \leq t\}$. Define, for some constants $\delta, K$ to be chosen later

$$(2.35) \qquad T(x^0, \mathcal{R}) = \inf\left\{t : \|\mathcal{K}(S^f(t))\| \leq \frac{\|x^0\|(1+\delta)^{n_t}}{K}\right\},$$

where $\mathcal{K}(S^f(0)) = x^0$. First, our goal is to prove that $\mathbb{E}(T(x^0, \mathcal{R}))^s = \infty$ if $K$ is large and $\delta$ is small enough.

From (2.9), Condition 1.1 and Lemma 2.3, we obtain that, for some $C_1$

$$(2.36) \quad T(x^0, \mathcal{R}) \geq C_1 \|(\tilde{A}(0) + \tilde{A}(1)\tilde{A}(0) + \cdots + \tilde{A}(U)\ldots\tilde{A}(0))x^0\|,$$

where

$$(2.37) \qquad U = \min\left\{k : \left\|\frac{\tilde{A}(k)}{(1+\delta)} \cdots \frac{\tilde{A}(0)}{(1+\delta)} x^0\right\| \leq \frac{\|x^0\|}{K}\right\}.$$

Denote for any $a > 0$

$$T^a(\mathcal{R}) = \inf_{x:\|x\|=a} T(x, \mathcal{R}), \qquad \hat{T}^a(\mathcal{R}) = \sup_{x:\|x\|=a} T(x, \mathcal{R}).$$

Using Condition 1.1 and Lemma 2.3, we obtain that there exists $C_2 > 0$ such that

$$(2.38) \qquad \frac{\hat{T}^a(\mathcal{R})}{T^a(\mathcal{R})} \leq C_2$$

for any $a$ and $\mathcal{R}$. By (2.36) and Lemma 2.3,

$$T(x^0, \mathcal{R}) \geq C_1 \left\|\left(\frac{\tilde{A}(0)}{(1+\delta)} + \frac{\tilde{A}(1)\tilde{A}(0)}{(1+\delta)^2} + \cdots + \frac{\tilde{A}(U)\ldots\tilde{A}(0)}{(1+\delta)^U}\right)x^0\right\|$$

$$\geq C_3 \left\|\frac{\tilde{A}(0)}{(1+\delta)} + \frac{\tilde{A}(1)\tilde{A}(0)}{(1+\delta)^2} + \cdots + \frac{\tilde{A}(U)\ldots\tilde{A}(0)}{(1+\delta)^U}\right\| \|x^0\|$$

for some $C_3$.

If $K$ is large enough, then there exist $1 < K' \leq K$ and $C_4$ such that

$$(2.39) \qquad \begin{aligned} T^a(\mathcal{R}) &\geq C_4 a \left\|\frac{\hat{A}(0)}{(1+\delta)^{d+1}} + \frac{\hat{A}(1)\hat{A}(0)}{(1+\delta)^{2(d+1)}} + \cdots + \frac{\hat{A}(U')\ldots\hat{A}(0)}{(1+\delta)^{U'(d+1)}}\right\| \\ &=: \Psi^a, \end{aligned}$$

where [compare with (2.37)]

$$U' = \min\left\{k : \left\|\frac{\hat{A}(k)\ldots\hat{A}(0)}{(1+\delta)^{k(d+1)}}\right\| \leq \frac{1}{K'}\right\}. \tag{2.40}$$

Abbreviate $b = \|x^0\|$. Using (2.38) and (2.40), one can write

$$\left\|\frac{\hat{A}(0)}{(1+\delta)^{d+1}} + \frac{\hat{A}(1)\hat{A}(0)}{(1+\delta)^{2(d+1)}} + \frac{\hat{A}(2)\hat{A}(1)\hat{A}(0)}{(1+\delta)^{3(d+1)}}\cdots\right\|$$

$$\leq \left\|\frac{\hat{A}(0)}{(1+\delta)^{d+1}} + \cdots + \frac{\hat{A}(U')\ldots\hat{A}(0)}{(1+\delta)^{U'(d+1)}}\right\|$$

$$+ \frac{1}{K'}\left\|\frac{\hat{A}(U'+1)\hat{A}(U')\ldots\hat{A}(0)}{(1+\delta)^{(U'+1)(d+1)}} + \cdots\right\|.$$

Continuing in this way and using (2.39), we can write

$$\left\|\frac{\hat{A}(0)}{(1+\delta)^{d+1}} + \frac{\hat{A}(1)\hat{A}(0)}{(1+\delta)^{2(d+1)}} + \frac{\hat{A}(2)\hat{A}(1)\hat{A}(0)}{(1+\delta)^{3(d+1)}}\cdots\right\| \leq C_5 \sum_{n=1}^{\infty} \frac{Z_n}{(K')^n},$$

where $Z_i$, $i = 1, 2, 3, \ldots$, are independent copies of $\Psi^b$. When $\mathsf{k}(s) > 1$ and $\delta$ is small enough to ensure that $(1+\delta)^{-1}\mathsf{k}(s) > 1$, Lemma 2.1 implies that [use $\Theta_k = (1+\delta)^{-1}\hat{A}(k)$, so that $\mathsf{k}^{\Theta}(s) = (1+\delta)^{-1}\mathsf{k}(s)$]

$$\mathbb{E}\left\|\frac{\hat{A}(0)}{(1+\delta)^{d+1}} + \frac{\hat{A}(1)\hat{A}(0)}{(1+\delta)^{2(d+1)}} + \frac{\hat{A}(2)\hat{A}(1)\hat{A}(0)}{(1+\delta)^{3(d+1)}}\cdots\right\|^s = \infty,$$

and so, by Lemma 2.2, $\mathbb{E}\Psi^b = \infty$. From (2.39), we obtain that

$$\mathbb{E}(T(x^0, \mathcal{R}))^s = \infty. \tag{2.41}$$

Let [recall (2.35)] $\kappa = n_{T(x^0, \mathcal{R})}$. Denote

$$B_n = \left\{\|\mathcal{K}(S(\sigma_n))\| \geq \left(1 - \frac{\delta}{2}\right)\|\mathcal{K}(S^f(\sigma_n^f))\|\right\}$$

and

$$B = \bigcap_{n=1}^{\kappa} B_n.$$

From part (ii) of Lemma 2.4 we obtain that for $n \leq \kappa$

$$\mathsf{P}_{\mathcal{R}}[B_n \mid B_{n-1}, \ldots, B_1] \geq 1 - \exp\{-C_6\|\mathcal{K}(S^f(\sigma_n^f))\|\} \tag{2.42}$$

for some $C_6 > 0$. By the definition (2.35) of $T(x^0, \mathcal{R})$,

$$\|\mathcal{K}(S^f(\sigma_n^f))\| \geq \frac{\|x^0\|(1+\delta)^n}{K},$$

so (2.42) implies that

$$P_{\mathcal{R}}[B] \geq \prod_{n=1}^{\infty}(1 - e^{-C_7(1+\delta)^n}) > 0.$$

Now, on $B$ we have that $\tau > T(x^0, \mathcal{R})$, so $\mathbb{E}_{\mathcal{R}}\tau^s \geq C_8(T(x^0, \mathcal{R}))^s$ for some positive constant $C_8$. Using (2.41), we obtain that

$$\mathbb{E}\tau^s \geq C_8 \mathbb{E}(T(x^0, \mathcal{R}))^s = \infty.$$

I. MacPhee  
M. Menshikov  
Mathematics Department  
University of Durham  
South Road  
Durham DH1 3LE  
United Kingdom  
E-mail: i.m.macphee@durham.ac.uk  
       mikhail.menshikov@durham.ac.uk

D. Petritis  
Institut de recherche mathématique  
Université de Rennes 1 and CNRS UMR 6625  
Campus de Beaulieu  
35042 Rennes Cedex  
France  
E-mail: dimitri.petritis@univ-rennes1.fr

S. Popov  
Instituto de matemática e estatística  
Universidade de São Paulo  
rua do Matão 1010  
CEP 05508-090  
São Paulo SP  
Brasil  
E-mail: popov@ime.usp.br